\documentclass[11pt]{amsart}

\usepackage{amsmath}
\usepackage{amssymb}
\usepackage{amsthm}
\usepackage{enumerate}
\usepackage{amsbsy}
\usepackage{amsfonts}

\usepackage[bitstream-charter]{mathdesign}
\usepackage[T1]{fontenc}

\usepackage{amsthm}
\newtheorem{theorem}{Theorem}
\newtheorem{proposition}{Proposition}
\newtheorem{lemma}{Lemma}
\newtheorem{corollary}{Corollary}
\theoremstyle{definition}

\theoremstyle{remark}
\newtheorem{remark}{Remark}

\newcommand{\bC}{\mathbb{C}}
\newcommand{\bE}{\mathbb{E}}
\newcommand{\bN}{\mathbb{N}}
\newcommand{\bP}{\mathbb{P}}
\newcommand{\bR}{\mathbb{R}}
\newcommand{\bS}{\mathbb{S}}
\newcommand{\bZ}{\mathbb{Z}}
\newcommand{\cE}{\mathcal{E}}
\newcommand{\cF}{\mathcal{F}}
\newcommand{\cM}{\mathcal{M}}

\newcommand{\De}{\Delta}
\newcommand{\de}{\delta}
\newcommand{\Gam}{\Gamma}
\newcommand{\gam}{\gamma}

\newcommand{\lam}{\lambda}

\newcommand{\si}{\sigma}
\newcommand{\s}{\sqrt}
\newcommand{\f}{\frac}
\newcommand{\set}[1]{{\{#1\}}}
\newcommand{\Set}[1]{{\left\{#1\right\}}}

\newcommand{\Abs}[1]{\Big\vert#1\Big\vert}
\newcommand{\Rd}{{\bR^d}}

\newcommand{\intzi}{\int_0^\infty}
\newcommand{\intzo}{\int_0^1}
\newcommand{\intii}{\int_{-\infty}^\infty}
\newcommand{\intoi}{\int_1^\infty}
\newcommand{\intRd}{\int_{\mathbb{R}^d}}
\newcommand{\snzi}{\sum_{n=0}^{\infty}}
\newcommand{\snoi}{\sum_{n=1}^{\infty}}

\begin{document}
\baselineskip=18pt

\title{Precise Asymptotic Approximations for Kernels Corresponding to L\'evy Processes}

\author{Sihun Jo \& Minsuk Yang}

\address{Sihun Jo: School of Mathematics, Korea Institute for Advanced Study, 85 Hoegiro Dongdaemungu, Seoul, Republic of Korea}
\address{Minsuk Yang: School of Mathematics, Korea Institute for Advanced Study, 85 Hoegiro Dongdaemungu, Seoul, Republic of Korea}
\email{yangm@kias.re.kr}

\begin{abstract}
By using basic complex analysis techniques, we obtain precise asymptotic approximations for kernels corresponding to symmetric $\alpha$-stable processes and their fractional derivatives.
We use the deep connection between the decay of kernels and singularities of the Mellin transforms.
The key point of the method is to transform the multi-dimensional integral to the contour integral representation.
We then express the integrand as a combination of gamma functions so that we can easily find all poles of the integrand.
We obtain various asymtotics of the kernels by using Cauchy's Residue Theorem with shifting contour integration. 
As a byproduct, exact coefficients are also obtained. 
We apply this method to general L\'evy processes whose characteristic functions are radial and satisfy some regularity and size conditions.
Our approach is based on the Fourier analytic point of view.
\\

\end{abstract}

\maketitle

\section{Introduction}
\label{S1}

This paper is concerned with the asymptotic behavior of the transition probability density of the L\'evy process in $\Rd$, $d\ge2$ with a radial L\'evy symbol.
The characteristic function of a L\'evy process $\set{X_t: t\ge0}$ has a specific form 
\[\bE(\exp(i\xi\cdot X_t)) = \exp(-t\eta(\xi))\]
where $\eta$ is a L\'evy symbol.
The symbol $\eta$ plays a central role.
Actually, it is possible to characterize all L\'evy processes by their L\'evy symbol $\eta$, in the sense that, one can construct a L\'evy process from any symbol of the form given in the L\'evy--Khinchine representation formula.
In Blumenthal-Getoor \cite{MR0123362} the index of $\eta$ given by
\[\alpha = \inf\Set{\lam>0 : \lim_{|\xi|\to\infty} \f{\eta(\xi)}{|\xi|^\lam}=0}\]
should satisfy $0 \le \alpha \le 2$.

Let us mention some of the fundamental symbols which are radial.
Brownian motion is a basic example which lies at the source of the theory of stochastic processes.
A symmetric $\alpha$-stable process is another important example whose symbol is given by
\[\eta(\xi)=|\xi|^\alpha, \quad 0<\alpha<2.\]
The infinitesimal generator of a symmetric $\alpha$-stable process $X_t$ is the fractional Laplacian $-(-\De)^{\alpha/2}$ that can be written in the following principal value integral
\[(-\De)^{\alpha/2}u(x) = C\mathrm{p.v.}\intRd \f{u(x)-u(y)}{|x-y|^{d+\alpha}} dy\]
for some constant $C$.
Recently there has been many interests from the theory of probability, PDE and SPDE to study such fractional Laplacians (see e.g., Caffarelli-Silvestre \cite{MR2354493}, Chen-Kim-Song \cite{MR2677618}, Bogdan-Sztonyk \cite{MR2320691}).
One more example is a relativistic stable process whose symbol is given by
\[\eta(\xi)=(|\xi|^2+m^2)^{\alpha/2} - m^\alpha.\]
It has no scaling property when $m\neq0$, but it coincides to the symmetric $\alpha$-stable process when $m=0$.
Sharp kernel estimates for relativistic stable processes in open sets are obtained by Chen-Kim-Song \cite{MR2917772}.

Analytic tools such as semigroup, resolvent and infinitesimal generator have deep connections with the probabilistic properties such as transience and recurrence, asymptotic behavior at infinity, etc.
The family of convolution operators indexed by $t\ge0$ and given by
\[P_tf(x) = \intRd f(x+y) \bP(X_t\in dy)\]
is a Markov semigroup.
Taking the Fourier transform, we write
\[\cF(P_tf)(\xi) = \bE \intRd e^{i\xi\cdot x} f(X_t+x) dx = e^{-t\eta(-\xi)} \cF f(\xi).\]
Thus, it is very important to obtain sharp estimates for kernels $P_t(x)$ given by
\[P_t(x) = (2\pi)^{-d} \intRd e^{-ix\cdot\xi} e^{-t\eta(-\xi)} d\xi.\]
For more information about L\'evy processes, see Sato's monograph  \cite{MR1739520}.

There are many results about asymptotic behavior of stable laws.
The asymptotic expansion of one-dimensional stable laws was obtained by Linnik \cite{MR0065055}, Skorohod \cite{MR0065839}, Zolotarev \cite{MR854867}, and others.
However, there are few results for the asymptotic expansion of stable laws in dimension $d$ greater than 1.
In \cite{MR1744782} Kolokoltsov studied the asymptotic behaviour of stable and stable-like processes and gave generalizations and applications.
Recently, derivatives of various fractional powers of the densities discussed and their expansions and bounds are also analyzed in detail in his book \cite{MR2780345}, and the main results of \cite{MR1744782} are now reproduced in a bit more cleaned-up and enlarged version in the Ch. 7 of that book. 

In this paper we propose a different method to obtain asymptotic approximations for $0<\alpha<2$
\[P_t^\alpha(x) = (2\pi)^{-d} \intRd e^{-ix\cdot\xi} e^{-t|\xi|^\alpha} d\xi\]
and their fractional derivatives for $\beta\ge0$
\[(-\De)^{\beta/2}P_t^\alpha(x) = (2\pi)^{-d} \intRd e^{-ix\cdot\xi} |\xi|^\beta e^{-t|\xi|^\alpha} d\xi.\]
Our argument of the proof is simple and direct, so we can apply the same method to general radial symbol $\eta$ which has some regularity conditions similar to the standard symbol class in the theory of pseudo-differential operators. As a byproduct, exact coefficients are also obtained.

Some of our main results about asymptotic expansions of $\alpha$-stable processes were originally due to Kolokoltsov \cite{MR1744782}.
We think that the underlying principle is the same, but its realization is different.
Our approach is based on the Fourier analytic point of view.
We try to use effectively the deep connection between the decay of kernels and singularities of the Mellin transforms.
The main tool is based on the very basic but decisive complex analysis technique, that is Cauchy's Residue Theorem with shifting contour integration.
The key point of the method is to transform the multi-dimensional integral to the contour integral representation.
We then express the integrand as a combination of gamma functions so that we can easily find all poles of the integrand.

\section{Preliminaries}
\label{S2}
In this section we introduce the Fourier transform of radial functions, the Mellin transform, and Mellin's inversion formula.
The gamma function is an important example of Mellin transform and the Bessel functions naturally arise from the spherical symmetry.
For reader's convenience, we give some of the basic properties of gamma and Bessel functions used in the proof.

\subsection{Fourier Transform}
\label{S2-1}
Given an integrable function $f(x)$ on $\Rd$, we define its Fourier transform $\cF f(\xi)$ by
\[\cF f(\xi) = \intRd e^{i\xi\cdot x} f(x) dx,\]
where $\xi\in\Rd$ and $\xi\cdot x=\sum_{j=1}^\infty \xi_j x_j$.
Similarly the inverse Fourier transform is defined by
\[\cF^{-1} f(x) = \f{1}{(2\pi)^d} \intRd e^{-ix\cdot\xi} f(\xi) d\xi.\]
Moreover, if $f$ and $\cF f$ are integrable on $\Rd$, then the Fourier inversion formula holds, that is,
\[f(x) = \cF^{-1} (\cF f)(x).\]

More generally, the space of integrable functions on $\Rd$ is contained in the space of finite complex measures on $\Rd$ with the total variation norm via the identification $fdx=d\mu$.
Thus the definition of Fourier transform is extended to measures via
\[\cF \mu(\xi)=\intRd e^{i\xi\cdot x} d\mu(x).\]
We now present a well-known lemma for the Fourier transform of radial functions.
We include the proof for reader's convenience.

\begin{lemma}
If $f$ is radial and integrable on $\Rd$, then
\begin{equation}\label{S2-1:E1}
\intRd e^{-ix\cdot\xi} f(|\xi|) d\xi = (2\pi)^{d/2} |x|^{-d/2+1} \intzi J_{d/2-1}(|x|r) f(r) r^{d/2} dr.
\end{equation}
\end{lemma}

\begin{proof}
In polar coordinates,
\begin{align*}
\intRd e^{-ix\cdot\xi} f(|\xi|) d\xi 
&= \intzi \left(\int_{\bS^{d-1}} e^{-ix\cdot r\theta} d\si(\theta)\right) f(r) r^{d-1} dr \\
&= \intzi \cF \si(-rx) f(r) r^{d-1} dr
\end{align*}
where $\si$ denotes surface measure on $\bS^{d-1}$.
The result follows from 
\begin{equation}\label{S2-1:E2}
\cF \si(\xi) = (2\pi)^{d/2} |\xi|^{-d/2+1} J_{d/2-1}(|\xi|),
\end{equation}
where $J_\nu$ is the Bessel function of the first kind of order $\nu$.
\qed
\end{proof}

\begin{remark}
For the proof of \eqref{S2-1:E2} and much more general information about Fourier transforms of measures supported on surfaces, we refer to Stein's monograph \cite{MR1232192}.
\end{remark}

\subsection{Mellin Transform}
\label{S2-2}
Given a function $g(r)$ on $[0,\infty)$, we define its Mellin transform $\cM g(z)$ by
\begin{equation}\label{S2-2:E1}
\cM g(z)=\intzi g(r) r^{z-1} dr.
\end{equation}
This integral defines a holomorphic function on a strip in the complex plane.
It is closely connected to the Fourier transform.
Actually, it can be expressed as a Fourier transform.
To see this, let $z=a+ib$ and then change variables $r=e^{-x}$ to obtain
\[\cM g(a+ib)=\intii g(e^{-x}) e^{-ax} e^{-ibx} dx.\]
From the Fourier inversion formula, we obtain Mellin's inversion formula
\begin{equation}\label{S2-2:E2}
g(r) = \f{1}{2\pi i} \int_{\Re(z)-i\infty}^{\Re(z)+i\infty} \cM g(z) r^{-z} dz.
\end{equation}

From now, we shall use the notation, for $c\in\bR$, 
\begin{equation}\label{S2-2:E3}
\int_{(c)} f(z) dz := \int_{c-i\infty}^{c+i\infty} f(z) dz.
\end{equation}

\subsection{Gamma Function}
\label{S2-3}
The gamma function is the Mellin transform of $e^{-t}$.
It is initially defined for $\Re(z)>0$ by the absolutely convergent integral
\begin{equation}\label{S2-3:E1}
\Gam(z) = \intzi e^{-t} t^{z-1} dt.
\end{equation}

Integrating by parts gives the functional equation 
\begin{equation}\label{S2-3:E2}
\Gam(z+1)=z\Gam(z)
\end{equation}
and this relation produces a meromorphic continuation.
We note that there is another method to obtain the meromorphic continuation.
Splitting the integral and expanding $e^{-t}$ in a power series produces
\begin{align*}
\Gam(z)
&= \intzo \snzi \f{(-t)^n}{n!} t^{z-1} dt + \intoi e^{-t} t^{z-1} dt \\
&= \snzi \f{(-1)^n}{n!(z+n)} + \intoi e^{-t} t^{z-1} dt.
\end{align*}
The series defines a meromorphic function with simple poles at $z=0,-1,-2,\dots$ and the last integral defines an entire function.
Thus, the gamma function has an analytic continuation to a meromorphic function on $\bC$ whose singularities are simple poles and the residue at $z=-n$ is $(-1)^n/n!$.

We recall the following basic identities for the gamma function
\begin{align}
&\Gamma(z)\Gamma(1-z)=\frac{\pi}{\sin(\pi z)},\label{S2-3:E3} \\
&\Gamma(z)\Gamma(-z)=-\frac{\pi}{z}\frac{1}{\sin(\pi
z)},\label{S2-3:E4} \\
&\Gamma(2z)=\frac{2^{2z-1}}{\sqrt{\pi}}\Gamma(z)\Gamma(z+\frac{1}{2}).\label{S2-3:E5} 
\end{align}

From Euler's reflection formula \eqref{S2-3:E3}, the function $1/\Gam(z)$ is entire with simple zeros at $z=0,-1,-2,\dots$ and it vanishes nowhere else.

We shall frequently use the following version of Stirling's formula for the gamma function:
\begin{equation}\label{S2-3:E6}
\Gamma(z)=e^{-z}e^{(z-\frac{1}{2})\log z}(2\pi)^{\frac{1}{2}}
\left(1+\frac{1}{12z}+\frac{1}{288z^2}-\frac{139}{51840z^3}+O(\frac{1}{z^4})\right)
\end{equation}
if $|\arg z|<\pi$.
In particular, 
\begin{equation}\label{S2-3:E7}
\lim_{|y|\rightarrow
\infty}|\Gamma(x+iy)|e^{\frac{1}{2}\pi|y|}|y|^{\frac{1}{2}-x}=(2\pi)^{\frac{1}{2}}.
\end{equation}

\subsection{Bessel Function}
\label{S2-4}
The Bessel function of the first kind of order $\nu$ can be defined by a number of different ways.
One of the usual ways is the Poisson representation formula which is given by
\[J_\nu(r) = \f{(r/2)^\nu}{\Gam(\nu+1/2)\Gam(1/2)} \int_{-1}^1 e^{irs} (1-s^2)^{\nu-1/2} ds,\]
where $\Re(\nu)>-1/2$ and $r\ge0$.
We note that
\begin{equation}\label{S2-4:E1}
|J_\nu(r)| \lesssim \min\set{r^\nu, r^{-1/2}}
\end{equation}
and 
\begin{equation}\label{S2-4:E2}
D[r^\nu J_\nu(r)] = r^\nu J_{\nu-1}(r)
\end{equation}
where $D$ denotes the derivative with respect to $r$.

We shall frequently use the Mellin transform of $r^{-\nu}J_{\nu}(r)$, that is,
\begin{equation}\label{S2-4:E3}
\intzi r^{-\nu}J_{\nu}(r) r^{z-1} dr = \f{2^{z-\nu-1}\Gam(\f{1}{2}z)}{\Gam(\nu-\f{1}{2}z+1)}
\end{equation}
for $\nu>-1/2$ and for $0<\Re(z)<\Re(\nu)+3/2$.
For more information about Bessel functions, see Watson's monograph \cite{MR1349110}.

In this paper we denote $f \lesssim g$ if $f \le cg$ for some constant $c>0$ and denote $f \approx g$ if $f \lesssim g$ and $g \lesssim f$.
We use the standard notation $f=O(g)$ and $f \sim g$.

\section{Symmetric Alpha-Stable Processes}
\label{S3}
In this section we prove complete aymptotics of kernels corresponding to symmetric $\alpha$-stable processes and their fractional derivatives via Mellin's trasform and give some corollaries.
The method of proof shows that the asymptotic behavior of an oscillatory integral is closely connected with the singularities of the Mellin transform of that integral.

\subsection{Asymptotic Expansion for Symmetric Alpha-Stable Process}
\label{S3-1}
We begin by considering the classical heat kernel
\[P_t(x) = (2\pi)^{-d} \intRd e^{-ix\cdot\xi} e^{-t|\xi|^2} d\xi.\]
It has the scaling property
\[P_t(x) = t^{-d/2} P_1(t^{-1/2}x).\]
This reflects selfsimilarity of the Brownian motion, that is, any change of time scale for the Brownian motion has the same effect as some change of spatial scale.
Thus we may assume $t=1$ and use Fubini's theorem to obtain
\[P_1(x) = (2\pi)^{-d} \prod_{n=1}^d \intii e^{-ix_n\xi_n} e^{-\xi_n^2} d\xi_n.\]
The function $e^{-z^2}$ is holomorphic in the entire complex plane and has rapid decay as $|\Re(z)|\to\infty$.
Hence by Cauchy's theorem we can shift the contour integration to yield
\[\intii e^{-ix_n\xi_n} e^{-\xi_n^2} d\xi_n = e^{-x_n^2/4} \intii e^{-u^2} du = \s{\pi} e^{-x_n^2/4}.\]
Therefore 
\[P_t(x) = (4\pi t)^{-d/2} e^{-|x|^2/4t}.\]

There are two key ingredients in the above argument.
The first one is to transform the multi-dimensional integral to the contour integration.
The second key ingredient is to evaluate the integral by shifting the contour.
We shall use the above fundamental technique to general cases.

We now consider the transition density of the symmetric $\alpha$-stable processes
\begin{equation}\label{S3-1:E1}
P_t^\alpha(x) = (2\pi)^{-d} \intRd e^{-ix\cdot\xi} e^{-t|\xi|^\alpha} d\xi.
\end{equation}

\begin{proposition}\label{S3-1:T1}
Let $d\ge2$ and $0<\alpha\le2$.
For all $t>0$, $x\in\Rd$, and $(d-1)/2<c<d$,
\begin{equation}\label{S3-1:E2}
P_t^\alpha(x)=\f{t^{-d/\alpha}}{\alpha \pi^{d/2}} \frac{1}{2\pi i} \int_{(c)} \f{\Gam(\f{z}{\alpha}) \Gam(\f{d-z}{2}) 2^{-z}}{\Gam(\f{z}{2})} (t^{-1/\alpha}|x|)^{-d+z} dz.
\end{equation}
\end{proposition}

\begin{proof}
We may assume that $t=1$ because of the scaling property
\begin{equation}\label{S3-1:E3}
P_t^\alpha(x) = t^{-d/\alpha} P_1^\alpha(t^{-1/\alpha}x).
\end{equation}
Using \eqref{S2-1:E1}, we have
\begin{equation}\label{S3-1:E4}
P_1^\alpha(x) = (2\pi)^{-d/2} |x|^{-d/2+1} \intzi J_{d/2-1}(|x|r) e^{-r^\alpha} r^{d/2} dr.
\end{equation}

Since the gamma function is the Mellin transform of $e^{-r}$, Mellin's inversion formula \eqref{S2-2:E2} implies that for all $c>0$,
\[e^{-r} = \frac{1}{2\pi i} \int_{(c)} \Gam(z) r^{-z} dz.\]
By rescaling we obtain the key identity
\begin{equation}\label{S3-1:E5}
e^{-r^\alpha} = \frac{1}{2\pi i} \int_{(c)} \alpha^{-1} \Gam(\alpha^{-1}z) r^{-z} dz.
\end{equation}

Inserting \eqref{S3-1:E5} to \eqref{S3-1:E4} and scaling, we have
\begin{align*}
P_1^\alpha(x)
&= (2\pi)^{-d/2} |x|^{-d/2+1} \intzi J_{d/2-1}(|x|r) \frac{1}{2\pi i} \int_{(c)} \alpha^{-1} \Gam(\alpha^{-1}z) r^{-z} dz r^{d/2} dr \\
&= (2\pi)^{-d/2} |x|^{-d/2+1} \frac{1}{2\pi i} \int_{(c)} \alpha^{-1} \Gam(\alpha^{-1}z) \intzi J_{d/2-1}(|x|r) r^{-z+d/2} dr dz \\
&= (2\pi)^{-d/2} |x|^{-d} \frac{1}{2\pi i} \int_{(c)} \alpha^{-1} \Gam(\alpha^{-1}z) \intzi J_{d/2-1}(r) r^{-z+d/2} dr |x|^{z} dz.
\end{align*}
The identity \eqref{S2-4:E3} gives for $(d-1)/2<c<d$
\[\intzi J_{d/2-1}(r) r^{-z+d/2} dr = \f{2^{-z+d/2}\Gam(\f{d-z}{2})}{\Gam(\f{z}{2})},\]
and therefore
\begin{equation}\label{S3-1:E6}
P_1^\alpha(x) = \f{1}{\alpha \pi^{d/2}} \frac{1}{2\pi i} \int_{(c)} \f{\Gam(\f{z}{\alpha}) \Gam(\f{d-z}{2}) 2^{-z}}{\Gam(\f{z}{2})} |x|^{-d+z} dz.
\end{equation}
This completes the proof.
\qed
\end{proof}

We now derive full asymptotic expansions for small and large distance.
The following theorem is an immediate consequence of the contour integral representation \eqref{S3-1:E2}.
It is the same expansion derived in \cite{MR1744782}, (Proposition 2.1 and Proposition 2.2), with adjusting coefficients.
Our intention to reprove some of the propositions in \cite{MR1744782} is to illustrate our idea how to obtain various expansion from the representation formula \eqref{S3-1:E2}.
The argument is very simple and direct.
In order to justify the convergence issue, we shall give rather strict error bounds.
It is also straightforward but a little messy.
In the next section we will repeat essentially the same argument to extend the results to more general radial symbols.

\begin{theorem}\label{S3-1:T2}
Let $d\ge2$.
\begin{enumerate}
\item 
If $1<\alpha<2$, then the following asymptotic expansion
\begin{equation}\label{S3-1:E7}
P_t^\alpha(x)=\f{t^{-d/\alpha}}{\alpha \pi^{d/2}} \snoi \f{(-1)^n}{n!} \f{\Gam(\f{d+2n}{\alpha}) 2^{-d-2n}}{\Gam(\f{d+2n}{2})} (t^{-1/\alpha}|x|)^{2n}.
\end{equation}
converges absolutely for all $t>0$ and $x\in\Rd$.

If $\alpha=1$, then it converges absolutely for all $|x|/t<1$.
\item 
If $0<\alpha<1$, then the following asymptotic expansion
\begin{equation}\label{S3-1:E8}
P_t^\alpha(x)=\f{t^{-d/\alpha}}{\pi^{d/2}} \snoi (1-1_\bZ(n\alpha/2)) \f{(-1)^n}{n!} \f{\Gam(\f{d+n\alpha}{2}) 2^{n\alpha}}{\Gam(-\f{n\alpha}{2})} (t^{-1/\alpha}|x|)^{-d-n\alpha}.
\end{equation}
converges absolutely for all $t>0$ and $x\in\Rd$, or in $t^{-1/\alpha}|x|<1$.

If $\alpha=1$, then it converges absolutely for all $|x|/t>1$.
\end{enumerate}
\end{theorem}

\begin{proof}
We assume without loss of generality that $t=1$.
\begin{enumerate}
\item 
\textbf{[Step 1]}
We first recall the integral representation formula \eqref{S3-1:E2} with $t=1$
\[P_1^\alpha(x)=\f{1}{\alpha \pi^{d/2}} \frac{1}{2\pi i} \int_{(c)} \f{\Gam(\f{z}{\alpha}) \Gam(\f{d-z}{2}) 2^{-z}}{\Gam(\f{z}{2})} |x|^{-d+z} dz\]
where $(d-1)/2<c<d$.
We focus on the integrand
\[\f{\Gam(\f{z}{\alpha}) \Gam(\f{d-z}{2}) 2^{-z}}{\Gam(\f{z}{2})} |x|^{-d+z},\]
which is a meromorphic function.
We note that for $\Re(z)>(d-1)/2$, it has simple poles only at $z=d+2n$ where $n=0,1,2,\ldots$ because of the term $\Gam(\f{d-z}{2})$.

\textbf{[Step 2]}
Fix $\tilde{c}>d$ and define the contour $\gam_R$ consists of a rectangle with vertices $c-iR, c+iR, \tilde{c}+iR, \tilde{c}-iR$ with the clockwise orientation.
Then we have
\begin{align*}
\lim_{R\to\infty} \int_{c+iR}^{\tilde{c}+iR} \f{\Gam(\f{z}{\alpha}) \Gam(\f{d-z}{2}) 2^{-z}}{\Gam(\f{z}{2})} |x|^{-d+z} dz = 0, \\
\lim_{R\to\infty} \int_{\tilde{c}-iR}^{c-iR} \f{\Gam(\f{z}{\alpha}) \Gam(\f{d-z}{2}) 2^{-z}}{\Gam(\f{z}{2})} |x|^{-d+z} dz = 0,
\end{align*}
since the integrand has rapid decay as $|\Im(z)|\to\infty$ by Stirling's formula \eqref{S2-3:E7}.
Thus, using Cauchy's residue theorem, we can shift the vertical line integration to the right, $c \to \tilde{c}$ with calculating residues.
Hence we have the asymptotic representation
\begin{equation}\label{S3-1:E9}
P_1^\alpha(x)=\f{1}{\alpha \pi^{d/2}} \left(\sum_{n=0}^N \f{(-1)^n}{n!} \f{\Gam(\f{d+2n}{\alpha}) 2^{-d-2n}}{\Gam(\f{d+2n}{2})} |x|^{2n} + \frac{1}{2\pi i} \int_{(c)} \f{\Gam(\f{z}{\alpha}) \Gam(\f{d-z}{2}) 2^{-z}}{\Gam(\f{z}{2})} |x|^{-d+z} dz\right)
\end{equation}
where $d+2N<c<d+2N+2$.

\textbf{[Step 3]}
We shall give an error bound for the remaining integral in \eqref{S3-1:E9}.
Let $z=c+yi$.
We claim that
\begin{equation}\label{S3-1:E10}
\int_{(c)}\frac{\Gamma(\frac{z}{\alpha})\Gamma(\frac{d-z}{2})2^{-z}}{\Gamma(\frac{z}{2})} |x|^{-d+z}dz \lesssim c^{-\frac{c(\alpha-1)}{2\alpha}} c^{\frac{d+1}{2}}\alpha^{-\frac{c}{\alpha}} e^{\frac{c(\alpha-1)}{\alpha}} |x|^{-d+c}.
\end{equation}
The implicit constant does not depend on $c$.

To see this, we have from the identity \eqref{S2-3:E4}
\[\Gamma(\frac{d-z}{2})=-\frac{2\pi}{(d-z)\sin(\frac{\pi(d-z)}{2})} \frac{1}{\Gamma(\frac{z-d}{2})}\]
and write
\[\frac{\Gamma(\frac{z}{\alpha})\Gamma(\frac{d-z}{2})2^{-z}}{\Gamma(\frac{z}{2})}=-\frac{2^{-z+1}\pi}{(d-z)\sin(\frac{\pi(d-z)}{2})}\frac{\Gamma(\frac{z}{\alpha})}{\Gamma(\frac{z}{2})\Gamma(\frac{z-d}{2})}.\]
Stirling's formula \eqref{S2-3:E6} yields
\begin{align*}
\frac{\Gamma(\frac{z}{\alpha})}{\Gamma(\frac{z}{2})\Gamma(\frac{z-d}{2})}&
=\frac{1}{(2\pi)^{\frac{1}{2}}}
\frac{e^{(\frac{z}{\alpha}-\frac{1}{2})\log(\frac{z}{\alpha})}}
{e^{(\frac{z}{2}-\frac{1}{2})\log(\frac{z}{2})}
e^{(\frac{z-d}{2}-\frac{1}{2})\log(\frac{z-d}{2})}}
\frac{e^{-\frac{z}{\alpha}}}{e^{-\frac{z}{2}}e^{-\frac{z-d}{2}}}(1+O(\frac{1}{z}))
\end{align*}
and it is easy to see that
\[\frac{2^{-z}}{z-d}\frac{e^{-\frac{z}{\alpha}}}{e^{-\frac{z}{2}}e^{-\frac{z-d}{2}}} \lesssim e^{\frac{c(\alpha-1)}{\alpha}}2^{-c}.\]
Now we have
\begin{equation}\label{S3-1:E11}
\begin{split}
&\bigg|\frac{1}{\sin(\frac{\pi(d-z)}{2})} \frac{e^{(\frac{z}{\alpha}-\frac{1}{2})\log(\frac{z}{\alpha})}}
{e^{(\frac{z}{2}-\frac{1}{2})\log(\frac{z}{2})}
e^{(\frac{z-d}{2}-\frac{1}{2})\log(\frac{z-d}{2})}}\bigg|\\
&= \frac{\exp\left(\frac{1}{2}(\frac{c}{\alpha}-\frac{1}{2})\log((\frac{c}{\alpha})^2+\frac{y^2}{\alpha^2})\right)}
{\exp\left(\frac{1}{2}(\frac{c}{2}-\frac{1}{2})\log((\frac{c}{2})^2+\frac{y^2}{4})\right)
\exp\left(\frac{1}{2}(\frac{c-d}{2}-\frac{1}{2})\log((\frac{c-d}{2})^2+\frac{y^2}{4})\right)} \\
&\quad \times \frac{1}{\sin(\frac{\pi(d-z)}{2})} \frac{\exp(-\frac{y}{\alpha}\arctan{\frac{y}{c}})}{\exp(-\frac{y}{2}\arctan{\frac{y}{c}})
\exp(-\frac{y}{2}\arctan{\frac{y}{c-d}})}.
\end{split}
\end{equation}
For the first term on the RHS of \eqref{S3-1:E11}, we have
\begin{align*}
&\frac{\exp\left(\frac{1}{2}(\frac{c}{\alpha}-\frac{1}{2})\log((\frac{c}{\alpha})^2+\frac{y^2}{\alpha^2})\right)}
{\exp\left(\frac{1}{2}(\frac{c}{2}-\frac{1}{2})\log((\frac{c}{2})^2+\frac{y^2}{4})\right)
\exp\left(\frac{1}{2}(\frac{c-d}{2}-\frac{1}{2})\log((\frac{c-d}{2})^2+\frac{y^2}{4})\right)}\\
&\leq\frac{\exp\left(\frac{c}{2\alpha}\log(c^2+y^2)\right)\alpha^{-\frac{c}{\alpha}}
((\frac{c}{\alpha})^2+\frac{y^2}{\alpha^2})^{-\frac{1}{4}}}
{\exp\left(\frac{c}{2}\log((c-d)^2+y^2)\right)2^{-c}
((\frac{c}{2})^2+\frac{y^2}{4})^{-\frac{1}{4}}((\frac{c-d}{2})^2+\frac{y^2}{4})^{-\frac{d+1}{4}}}\\
&\lesssim \left(\frac{(c^2+y^2)^{\frac{1}{\alpha}}}{((c-d)^2+y^2)}\right)^{\frac{c}{2}}
2^c\alpha^{-\frac{c}{\alpha}}\left(\frac{c^2}{4}+\frac{y^2}{4}\right)^{\frac{d+1}{4}}\\
&\lesssim \left((c-d)^2+y^2\right)^{-\frac{c}{2}(1-\frac{1}{\alpha})}
\left(1+\frac{2cd-d^2}{(c-d)^2+y^2}\right)^{\frac{c}{2}}
2^c\alpha^{-\frac{c}{\alpha}}c^{\frac{d+1}{2}}(y^2+2)^{\frac{d+1}{4}}\\
&\lesssim c^{-\frac{c(\alpha-1)}{2\alpha}}2^c\alpha^{-\frac{c}{\alpha}}c^{\frac{d+1}{2}}(y^2+2)^{\frac{d+1}{4}}.
\end{align*}
For the second term on the RHS of \eqref{S3-1:E11}, we have for $\alpha<2$,
\begin{align*}
&\frac{1}{|\sin(\frac{\pi(d-z)}{2})|} \times \frac{1}{\exp(-\frac{y}{2}\arctan{\frac{y}{c-d}})} \times \frac{\exp(-\frac{y}{\alpha}\arctan{\frac{y}{c}})}{\exp(-\frac{y}{2}\arctan{\frac{y}{c}})} \\
&\lesssim \exp(-\pi|y|/2) \times \exp(\pi|y|/4) \times 1 \\
&\le \exp(-\pi|y|/4).
\end{align*}
Therefore, we obtain the following error bound
\begin{align*}
&\int_{(c)}\frac{\Gamma(\frac{z}{\alpha})\Gamma(\frac{d-z}{2})2^{-z}}{\Gamma(\frac{z}{2})} |x|^{-d+z}dz \\
&\lesssim c^{-\frac{c(\alpha-1)}{2\alpha}}2^c\alpha^{-\frac{c}{\alpha}}c^{\frac{d+1}{2}} e^{\frac{c(\alpha-1)}{\alpha}}2^{-c} |x|^{-d+c} \int_{-\infty}^{\infty}(y^2+2)^{\frac{d+1}{4}} \exp(-\pi|y|/4) dy
\end{align*}
and the result \eqref{S3-1:E10} follows.

\textbf{[Step 4]}
If $1<\alpha<2$, then the RHS of \eqref{S3-1:E10} goes to zero for all $x\in\Rd$ as $c\to\infty$.
If $\alpha=1$, then the RHS of \eqref{S3-1:E10} goes to zero for $|x|<1$ as $c\to\infty$.
If we let $c\to\infty$, then the series in \eqref{S3-1:E9} converges absolutely by using Stirling's formula \eqref{S2-3:E6} and the remaining integral in \eqref{S3-1:E9} goes to zero.
Therefore, we obtain the desired expansion \eqref{S3-1:E7}.

\item 
In this case, the proof is almost the same except for shifting the vertical line integration to the left.

\textbf{[Step 1]}
We focus on the integrand
\[\f{\Gam(\f{z}{\alpha}) \Gam(\f{d-z}{2}) 2^{-z}}{\Gam(\f{z}{2})} |x|^{-d+z}.\]
For $\Re(z)<d$, the function $\Gam(\f{z}{\alpha})$ has simple poles at $z=0,-\alpha,-2\alpha,\ldots$ 
and the function $1/\Gam(\f{z}{2})$ has simple zeros at $z=0,-2,-4,\ldots$.
We notice that if $n\alpha \in 2\bZ$, then the pole and the zero at that point is cancelled. 

\textbf{[Step 2]}
We shift the vertical line integration to the left with calculating residues.
Hence we have the asymptotic representation
\begin{equation}\label{S3-1:E12}
\begin{split}
P_1^\alpha(x)
&=\f{1}{\pi^{d/2}} \sum_{n=0}^N (1-1_\bZ(n\alpha/2)) \f{(-1)^n}{n!} \f{\Gam(\f{d+n\alpha}{2}) 2^{n\alpha}}{\Gam(-\f{n\alpha}{2})} |x|^{-d-n\alpha} \\
&\quad+ \f{1}{\pi^{d/2}} \frac{1}{2\pi i} \int_{(c)} \f{\Gam(\f{z}{\alpha}) \Gam(\f{d-z}{2}) 2^{-z}}{\Gam(\f{z}{2})} |x|^{-d+z} dz
\end{split}
\end{equation}
where $-(N+1)\alpha<c<-N\alpha$.
If we put $c=-C$, then
\begin{align*}
&\int_{(c)}\frac{\Gamma(\frac{z}{\alpha})\Gamma(\frac{d-z}{2})2^{-z}}{\Gamma(\frac{z}{2})}
|x|^{-d+z}dz \\
&=\int_{(-C)}\frac{\Gamma(\frac{z}{\alpha})\Gamma(\frac{d-z}{2})2^{-z}}{\Gamma(\frac{z}{2})} |x|^{-d+z}dz \\
&=\int_{(C)}\frac{\Gamma(\frac{-z}{\alpha})\Gamma(\frac{d+z}{2})2^{z}}{\Gamma(\frac{-z}{2})} |x|^{-d+z}dz.
\end{align*}
In order to obtain the desired expansion \eqref{S3-1:E8}, we need to get an error bound of the last integral as $C\to\infty$.

\textbf{[Step 3]}
We claim that the last integral has the following error bound
\begin{equation}\label{S3-1:E13}
\int_{(C)}\frac{\Gamma(\frac{-z}{\alpha})\Gamma(\frac{d+z}{2})2^{z}}{\Gamma(\frac{-z}{2})} |x|^{-d-z}dz \lesssim C^{-\frac{C(1-\alpha)}{2\alpha}}C^{\frac{d-1}{2}}\alpha^{\frac{C}{\alpha}} e^{\frac{1-\alpha}{\alpha}C} |x|^{-d-C}.
\end{equation}
The implicit constant does not depend on $C$.

To see this, we have from the identity \eqref{S2-3:E4}
\begin{align*}
\Gamma(-\frac{z}{\alpha})
&=-\frac{\pi\alpha}{z\sin(\frac{\pi z}{\alpha})}\frac{1}{\Gamma(\frac{z}{\alpha})} \\
\Gamma(-\frac{z}{2})
&=-\frac{2\pi}{z\sin(\frac{\pi z}{2})}\frac{1}{\Gamma(\frac{z}{2})}
\end{align*}
and write
$$\frac{\Gamma(\frac{-z}{\alpha})\Gamma(\frac{d+z}{2})2^{z}}{\Gamma(\frac{-z}{2})}
=\frac{\Gamma(\frac{z}{2})\Gamma(\frac{d+z}{2})}{\Gamma(\frac{z}{\alpha})}
\frac{2^{z-1}\pi \sin(\frac{\pi z}{2})}{ \sin(\frac{\pi
z}{\alpha})}$$
Stirling's formula \eqref{S2-3:E6} yields
\[\frac{\Gamma(\frac{z}{2})\Gamma(\frac{z+d}{2})}{\Gamma(\frac{z}{\alpha})} 
= (2\pi)^{\frac{1}{2}} \frac{e^{(\frac{z}{2}-\frac{1}{2})\log(\frac{z}{2})}e^{(\frac{z+d}{2}-\frac{1}{2})\log(\frac{z+d}{2})}}{e^{(\frac{z}{\alpha}-\frac{1}{2})\log(\frac{z}{\alpha})} } \frac{e^{-\frac{z}{2}}e^{-\frac{z+d}{2}}}{e^{-\frac{z}{\alpha}}}(1+O(\frac{1}{z})),\]
and it is easy to see that
\[\frac{2^ze^{-\frac{z}{2}}e^{-\frac{z+d}{2}}}{e^{-\frac{z}{\alpha}}} \lesssim 2^C e^{\frac{1-\alpha}{\alpha}C}.\]
Let $z=C+yi$.
Now we have
\begin{equation}\label{S3-1:E14}
\begin{split}
&\bigg|\frac{\sin(\frac{\pi z}{2})}{\sin(\frac{\pi z}{\alpha})} \frac{e^{(\frac{z}{2}-\frac{1}{2})\log(\frac{z}{2})}e^{(\frac{z+d}{2}-\frac{1}{2})\log(\frac{z+d}{2})}}{e^{(\frac{z}{\alpha}-\frac{1}{2})\log(\frac{z}{\alpha})} }\bigg| \\
&= \frac{\exp\left(\frac{1}{2}(\frac{C}{2}-\frac{1}{2})\log((\frac{C}{2})^2+\frac{y^2}{4})\right) \exp\left(\frac{1}{2}(\frac{C+d}{2}-\frac{1}{2})\log((\frac{C+d}{2})^2+\frac{y^2}{4})\right)}{\exp\left(\frac{1}{2}(\frac{C}{\alpha}-\frac{1}{2})\log((\frac{C}{\alpha})^2+\frac{y^2}{\alpha^2})\right)} \\
&\quad \times \frac{\sin(\frac{\pi z}{2})}{\sin(\frac{\pi z}{\alpha})} \frac{\exp(-\frac{y}{2}\arctan{\frac{y}{C}})\exp(-\frac{y}{2}\arctan{\frac{y}{C+d}})}{\exp(-\frac{y}{\alpha}\arctan{\frac{y}{C}})}.
\end{split}
\end{equation}
For the first term on the RHS of \eqref{S3-1:E14}, we have
\begin{align*}
&\frac{\exp\left(\frac{1}{2}
(\frac{C}{2}-\frac{1}{2})\log((\frac{C}{2})^2+\frac{y^2}{4})\right)
\exp\left(\frac{1}{2}(\frac{C+d}{2}-\frac{1}{2})\log((\frac{C+d}{2})^2+\frac{y^2}{4})\right)}
{\exp\left(\frac{1}{2}(\frac{C}{\alpha}-\frac{1}{2})\log((\frac{C}{\alpha})^2+\frac{y^2}{\alpha^2})\right)}\\
&\leq\frac{\exp\left(\frac{C}{2}\log((C+d)^2+y^2)\right)2^{-C}
((\frac{C}{2})^2+\frac{y^2}{4})^{-\frac{1}{4}}((\frac{C+d}{2})^2+\frac{y^2}{4})^{\frac{d-1}{4}}}
{\exp\left(\frac{C}{2\alpha}\log(C^2+y^2)\right)\alpha^{-\frac{C}{\alpha}}
((\frac{C}{\alpha})^2+\frac{y^2}{\alpha^2})^{-\frac{1}{4}}}\\
&\lesssim
\left(\frac{((C+d)^2+y^2)}{(C^2+y^2)^{\frac{1}{\alpha}}}\right)^{\frac{C}{2}}
2^{-C}\alpha^{\frac{C}{\alpha}}\left(\frac{(C+d)^2}{4}+\frac{y^2}{4}\right)^{\frac{d-1}{4}}\\
&\lesssim \left(C^2+y^2\right)^{-\frac{C}{2}(\frac{1}{\alpha}-1)}
\left(1+\frac{2Cd+d^2}{C^2+y^2}\right)^{\frac{C}{2}}
2^{-C}\alpha^{\frac{C}{\alpha}}\left(\frac{(C+d)^2}{4}+\frac{y^2}{4}\right)^{\frac{d-1}{4}}\\
&\lesssim C^{-\frac{C(1-\alpha)}{2\alpha}}2^{-C}\alpha^{\frac{C}{\alpha}}
C^{\frac{d-1}{2}}(y^2+2)^{\frac{d-1}{4}}
\end{align*}
For the second term on the RHS of \eqref{S3-1:E14}, we have for $\alpha\le1$,
\begin{align*}
&\Abs{\frac{\sin(\frac{\pi z}{2})}{\sin(\frac{\pi z}{\alpha})}}
\times \frac{\exp(-\frac{y}{2}\arctan{\frac{y}{C}})}{\exp(-\frac{y}{\alpha}\arctan{\frac{y}{C}})} 
\times \exp(-\frac{y}{2}\arctan{\frac{y}{C+d}}) \\
&\lesssim \exp(-(1/\alpha-1/2)\pi|y|) \times \exp((1/\alpha-1/2)\pi|y|/2) \times 1 \\
&\le \exp(-\pi|y|/4).
\end{align*}
Therefore, we obtain the following error bound
\begin{align*}
&\int_{(C)}\frac{\Gamma(\frac{-z}{\alpha})\Gamma(\frac{d+z}{2})2^{z}}{\Gamma(\frac{-z}{2})} |x|^{-d-z}dz \\
&\lesssim C^{-\frac{C(1-\alpha)}{2\alpha}}2^{-C}\alpha^{\frac{C}{\alpha}} C^{\frac{d-1}{2}}2^C e^{\frac{1-\alpha}{\alpha}C} |x|^{-d-C}\int_{-\infty}^{\infty}(y^2+2)^{\frac{d-1}{4}} \exp(-\pi|y|/4) dy.
\end{align*}

\textbf{[Step 4]}
If $0<\alpha<1$, then the RHS of \eqref{S3-1:E13} goes to zero for all $x\in\Rd$ as $C\to\infty$.
If $\alpha=1$, then the RHS of \eqref{S3-1:E13} goes to zero for $|x|>1$ as $C\to\infty$.
If we let $C\to\infty$, then the series in \eqref{S3-1:E12} converges absolutely by using Stirling's formula \eqref{S2-3:E6} and the remaining integral in \eqref{S3-1:E12} goes to zero.
Therefore, we obtain the desired expansion \eqref{S3-1:E8}.
\end{enumerate}
This completes the proof.
\qed
\end{proof}

\begin{remark}
For the case $\alpha=1$ the asymptotic expansion coincides with Poisson's kernel.
In fact, we have from the previous theorem
\[P_1^1(x)=\frac{1}{\pi^{d/2}|x|^{d+1}} \snzi \frac{(-1)^{2n+1}\Gamma(n+\f{d+1}{2}) 2^{2n+1}}{(2n+1)!\Gamma(-n-\f{1}{2})} |x|^{-2n}.\]
A calculation gives
\[(2n+1)!\Gamma(-n-1/2) = 2^{2n+1} n! (-1)^{n+1} \Gam(1/2)\]
and so we have
\begin{align*}
P_1^1(x)
&=\frac{1}{\pi^{\f{d+1}{2}}|x|^{d+1}} \snzi \frac{(-1)^n\Gamma(n+\f{d+1}{2})}{n!} |x|^{-2n} \\
&=\frac{\Gamma(\f{d+1}{2})}{\pi^{\f{d+1}{2}}|x|^{d+1}} \snzi \binom{-\f{d+1}{2}}{n} |x|^{-2n} \\
&=\frac{\Gamma(\f{d+1}{2})}{\pi^{\f{d+1}{2}}} (1+|x|^2)^{-\f{d+1}{2}}.
\end{align*}
\end{remark}

\begin{remark}
For the exceptional case $\alpha=2$ the integral representation \eqref{S3-1:E2} is still true.
So we have
\[P_1^2(x) = \f{1}{2\pi^{d/2} |x|^d} \frac{1}{2\pi i} \int_{(c)} \Gam(\f{d-z}{2}) 2^{-z} |x|^{z} dz.\]
By Cauchy's residue theorem we can shift the contour integration to the right to obtain
\begin{align*}
P_1^2(x)
&= \f{1}{\pi^{d/2} |x|^d} \snzi \f{(-1)^n}{n!} (|x|/2)^{d+2n} \\
&= (4\pi)^{-d/2} e^{-|x|^2/4}.
\end{align*}
\end{remark}

\subsection{Asymptotic Expansion for Fractional Derivatives}
\label{S3-2}
The fractional derivatives of $P_t^\alpha(x)$ is defined at least formally by
\[(-\De)^{\beta/2}P_t^\alpha(x) = (2\pi)^{-d} \intRd e^{-ix\cdot\xi} |\xi|^\beta e^{-t|\xi|^\alpha} d\xi.\]
Sometimes it is necessary to estimate the derivatives of transition densities especially when one studies a certain type of stochastic partial differential equations.
For example, Chang-Lee \cite{MR2946443} obtained Sobolev and Besov estimates for an SPDE with a fractional Lapalacian operator from the decay estimates of $P_t^\alpha(x)$ and it's derivatives.

\begin{proposition}\label{S3-2:T1}
Let $d\ge2$ and $0<\alpha\le2$.
For all $t>0$, $x\in\Rd$, and $(d-1)/2<c<d$,
\begin{equation}\label{S3-2:E1}
(-\De)^{\beta/2}P_t^\alpha(x) = \f{t^{-(d+\beta)/\alpha}}{\alpha \pi^{d/2}} \frac{1}{2\pi i} \int_{(c)} \f{\Gam(\f{z}{\alpha}) \Gam(\f{d+\beta-z}{2}) 2^{-z}}{\Gam(\f{z-\beta}{2})} (t^{-1/\alpha}|x|)^{-d-\beta+z} dz.
\end{equation}
\end{proposition}

\begin{proof}
We may assume that $t=1$ because of the scaling property
\[(-\De)^{\beta/2}P_t^\alpha(x) = t^{-(d+\beta)/\alpha} (-\De)^{\beta/2}P_1^\alpha(t^{-1/\alpha}x).\]
Using \eqref{S2-1:E1}, we have
\[(-\De)^{\beta/2}P_1^\alpha(x) = (2\pi)^{-d/2} |x|^{-d/2+1} \intzi J_{d/2-1}(|x|r) e^{-r^\alpha} r^{d/2+\beta} dr.\]
By the same way as in the proof of Proposition \ref{S3-1:T1}, we have
\[(-\De)^{\beta/2}P_1^\alpha(x) = \f{1}{\alpha \pi^{d/2}} \frac{1}{2\pi i} \int_{(c)} \f{\Gam(\f{z}{\alpha}) \Gam(\f{d+\beta-z}{2}) 2^{-z+\beta}}{\Gam(\f{z-\beta}{2})} |x|^{-d-\beta+z} dz.\]
This completes the proof. \qed
\end{proof}

By exactly the same way as in the proof of Theorem \ref{S3-1:T2}, one can immediately derive various asymptotic expansions for the fractional derivatives of transition probability densities.
However, we only present the first order approximation for convenience.

\begin{theorem}\label{S3-2:T2}
Let $d\ge2$, $0<\alpha<2$ and $\beta\ge0$.
\begin{enumerate}
\item 
If $\beta\neq 0,2,4,\ldots$, then as $|x|\to\infty$ the first order approximation becomes
\[(-\De)^{\beta/2}P_1^\alpha(x) \sim -\f{2^{\beta-2}}{\pi^{d/2+2}\alpha} \frac{\Gam(\f{d+\beta}{2})}{\Gam(-\f{\beta}{2})} |x|^{-d-\beta}.\]
\item 
If $\beta=0,2,4,\ldots$, then as $|x|\to\infty$ the first order approximation becomes
\[(-\De)^{\beta/2}P_1^\alpha(x) \sim -\f{2^{\beta+\alpha-2}}{\pi^{d/2+2}\alpha} \frac{\Gam((d+\beta+\alpha)/2)}{\Gam(-(\beta+\alpha)/2)} |x|^{-d-\beta-\alpha}.\]
\end{enumerate}
\end{theorem}

\begin{proof}
The first order approximation follows from the first pole of the integrand in \eqref{S3-2:E1} on the negative real axis.
If $\beta$ is not an even integer, then $0$ is the first pole of the integrand.
If $\beta$ is an even integer, then $-\alpha$ is the first pole of the integrand.
We shift the line of integral to the left to get the result.
\qed
\end{proof}

\begin{remark}
Kim-Kim \cite{MR2869738} obtained an upper bound
\[|(-\De)^{\beta/2}P_1^\alpha(x)| \lesssim |x|^{-(d+\beta)}.\]
This was crucial for them to prove their generalization of a Littlwood-Paley operator for the fractional Laplacian.
From the above theorem we have
\[\lim_{|x|\to\infty} |x|^{d+\beta} |(-\De)^{\beta/2}P_1^\alpha(x)| \neq 0.\]
\end{remark}

\begin{corollary}\label{S3-2:T3}
Let $d\ge2$, $0<\alpha<2$ and $\beta\ge0$.
If $\beta \neq 0,2,4,\dots$, then
\[|(-\De)^{\beta/2}P_t^\alpha(x)| \lesssim \min\Set{t^{-(d+\beta)/\alpha}, \f{1}{|x|^{d+\beta}}}.\]
If $\beta = 0,2,4,\dots$, then
\[|(-\De)^{\beta/2}P_t^\alpha(x)| \lesssim \min\Set{t^{-(d+\beta)/\alpha}, \f{t}{|x|^{d+\beta+\alpha}}}.\]
\end{corollary}

\begin{proof}
Since $|\xi|^\beta e^{-|\xi|^\alpha}$ is integrable on $\Rd$ for all $\beta\ge0$, the function $(-\De)^{\beta/2} p(1,x)$ is bounded and smooth on $(0,\infty)\times\Rd$ by the Riemann-Lebesgue lemma.
\qed
\end{proof}

\subsection{Sum of Two Powers}
\label{S3-3}
In this subsection we consider the following symbol
\[\eta(\xi) = |\xi|^\alpha+|\xi|^\beta\]
with $0<\alpha<\beta<2$.

A serious investigation of transition densities for symbols that are sums of two powers was initiated by Kochubei \cite{MR972089}.
Recently, there are many results in this subject.
The two-sided sharp estimates on the transition density function of the independent sum of two stable processes is given in \cite{MR2357678}.
Their upper bound estimates in particular imply the following Corollary.
We present the following simple computation because it gives us an initial hint, that is, the intensity of the singularity of the symbol at the origin is related to the decay of the kernel.
We note that the function $|\xi|^\alpha$ is smooth away from the origin and is H\"older continuous at the origin.

\begin{proposition}\label{sec3-thm3}
Let $d\ge2$, $0<\alpha<\beta<2$ and
\[K_t(x) = (2\pi)^{-d} \intRd e^{-ix\cdot\xi} \exp(-t[|\xi|^\alpha+|\xi|^\beta]) d\xi.\]
Then
\[\begin{cases}
K_t(x) \lesssim t^{-d/\beta} (1+t^{-1/\beta}|x|)^{-(d+\alpha)} & \text{ if } t\le1 \\
K_t(x) \lesssim t^{-d/\alpha} (1+t^{-1/\alpha}|x|)^{-(d+\alpha)} & \text{ if } t\ge1. \\
\end{cases}\]
\end{proposition}

\begin{proof}
Note that
\[K_t(x) = \intRd P_t^\alpha(x-y) P_t^\beta(y) dy.\]
It is easy to see that the first order approximation becomes
\begin{equation}\label{S3-2:E2}
P_t^\alpha(x) \approx t^{-d/\alpha} (1+t^{-1/\alpha}|x|)^{-(d+\alpha)}.
\end{equation}

First we assume $t\le1$ and consider the cases $t^{-1/\beta}|x|\le1$ and $t^{-1/\beta}|x|\ge1$.
In the case $t^{-1/\beta}|x|\le1$ we have
\[P_t^\beta(y) \lesssim t^{-d/\beta} \le \f{2^{d+\alpha} t^{-d/\beta}}{(1+t^{-1/\beta}|x|)^{d+\alpha}}.\]
Therefore
\[K_t(x) \lesssim \f{2^{d+\alpha} t^{-d/\beta}}{(1+t^{-1/\beta}|x|)^{d+\alpha}} \intRd  P_t^\alpha(x-y) dy.\]
In the case $t^{-1/\beta}|x|\ge1$ we split the integral to obtain
\begin{align*}
K_t(x)
&\lesssim \int_{|y-x|\ge|x|/2} P_t^\alpha(x-y) P_t^\beta(y) dy \\
&\quad + \int_{|y|\ge|x|/2} P_t^\alpha(x-y) P_t^\beta(y) dy \\
&\le \f{2^{d+\alpha} t^{-d/\alpha}}{(1+t^{-1/\alpha}|x|)^{d+\alpha}} \int_{|y-x|\ge|x|/2} P_t^\beta(y) dy \\
&\quad + \f{2^{d+\beta} t^{-d/\beta}}{(1+t^{-1/\beta}|x|)^{d+\beta}} \int_{|y|\ge|x|/2} P_t^\alpha(x-y) dy.
\end{align*}
Since $t^{-1/\beta}|x|\ge1$ and $t^{1-\alpha/\beta}\le1$, we have
\[\f{2^{d+\alpha}t^{-d/\alpha}}{(1+t^{-1/\alpha}|x|)^{d+\alpha}} \le \f{2^{d+\alpha}t}{|x|^{d+\alpha}} = \f{2^{d+\alpha}t^{-d/\beta+1-\alpha/\beta}}{(t^{-1/\beta}|x|)^{d+\alpha}} \le \f{4^{d+\alpha}t^{-d/\beta}}{(1+t^{-1/\beta}|x|)^{d+\alpha}}\]
and the required estimate follows.

Now we assume $t\ge1$ and consider the cases $t^{-1/\alpha}|x|\le1$ and $t^{-1/\alpha}|x|\ge1$.
In the case $t^{-1/\alpha}|x|\le1$ we have
\[P_t^\alpha(x-y) \le t^{-d/\alpha} \le \f{2^{d+\alpha} t^{-d/\alpha}}{(1+t^{-1/\alpha}|x|)^{d+\alpha}}\]
Therefore
\[K_t(x) \lesssim \f{2^{d+\alpha} t^{-d/\alpha}}{(1+t^{-1/\alpha}|x|)^{d+\alpha}} \intRd  P_t^\beta(y) dy.\]
In the case $t^{-1/\alpha}|x|\ge1$ we split the integral as before to obtain
\begin{align*}
K_t(x)
&\lesssim \f{2^{d+\alpha} t^{-d/\alpha}}{(1+t^{-1/\alpha}|x|)^{d+\alpha}} \int_{|y-x|\ge|x|/2} P_t^\beta(y) dy \\
&\quad + \f{2^{d+\beta} t^{-d/\beta}}{(1+t^{-1/\beta}|x|)^{d+\beta}} \int_{|y|\ge|x|/2} P_t^\alpha(x-y) dy.
\end{align*}
Since $t^{-1/\alpha}|x|\ge1$ and $t^{1-\beta/\alpha}\le1$, we have
\[\f{2^{d+\beta} t^{-d/\beta}}{(1+t^{-1/\beta}|x|)^{d+\beta}} \le \f{2^{d+\beta} t}{|x|^{d+\beta}} \le \f{4^{d+\beta} t^{-d/\alpha+1-\beta/\alpha}}{(1+t^{-1/\alpha}|x|)^{d+\beta}} \le \f{4^{d+\beta} t^{-d/\alpha}}{(1+t^{-1/\alpha}|x|)^{d+\alpha}}\]
and the required estimate follows.
\qed
\end{proof}

\begin{remark}
The estimate \eqref{S3-2:E2} in fact is due to Blumenthal and Getoor. 
The first term of the asymptotic behavior of $P^\alpha_t (x)$ is derived in Theorem 2.1 of \cite{MR0119247} from which the estimate \eqref{S3-2:E2} follows immediately.
\end{remark}

\begin{remark}
Local principles of approximate scaling and approximate independence of increments are studied in \cite{MR1744782}.
\end{remark}

\section{Generalizations to Radial Symbols}
\label{S4}
In this section by refining the same techniques already developed in the previous section, we investigate the asymptotic behavior of the kernels corresponding to more general radial symbol $\eta$ under some weak conditions.

\subsection{The First General Case}
\label{S4-1}
Let
\[K^\beta(t,x) = (2\pi)^{-d} \intRd e^{-ix\cdot\xi} |\xi|^\beta e^{-t\eta(|\xi|)} d\xi\]
for $\beta\ge0$.

To prove the following proposition, the key point is to find a Mellin's inversion formula for the symbol $e^{-t\eta(r)}$ and then to obtain a meromorphic continuation, which can be obtained by performing an integration by parts.

\begin{proposition}
Suppose that $\eta(r)$ satisfies
\begin{equation}\label{S4-1:E1}
\liminf_{r\to\infty} \f{\eta(r)}{\log r} = \infty,
\end{equation}
and for some $0<\alpha<2$ and $k>(d+3)/2+\beta$
\begin{equation}\label{S4-1:E2}
\sup_{0\le m\le k} \sup_r |r^{-\alpha+m} D^m \eta(r)| < \infty,
\end{equation}
where $D$ denotes the derivative with respect to $r$.
Then for all $t>0$, $x\in\Rd$ and $(d+1)/2+\beta<c<d+\beta$
\begin{equation}\label{S4-1:E3}
K^\beta(t,x) = \f{(-1)^k}{\pi^{d/2} |x|^{d+\beta}} \frac{1}{2\pi i} \int_{(c)} \f{\Gam(z)\Gam(\f{-z+d+\beta}{2})2^{-z+\beta}}{\Gam(z+k)\Gam(\f{z-\beta}{2})} M_t^k(z) |x|^z dz
\end{equation}
where
\[M_t^k(z) = \intzi D^k(e^{-t\eta(r)}) r^{z+k-1} dr.\]
\end{proposition}

\begin{proof}
\textbf{[Step 1]}
Let $M_t(z)$ denote the Mellin transform
\[M_t(z) = \intzi e^{-t\eta(r)} r^{z-1} dr.\]
This integral converges absolutely for $\Re(z)>0$ by the condition \eqref{S4-1:E1}.
Integrating by parts yields
\begin{align*}
M_t(z)
&= \intzi e^{-t\eta(r)} r^{z-1} dr \\
&= \f{-1}{z} \intzi D(e^{-t\eta(r)}) r^{z} dr \\
&= \f{(-1)^2}{(z+1)z} \intzi D^2(e^{-t\eta(r)}) r^{z+1} dr \\
&\qquad\qquad\qquad\qquad \vdots \\
&= \f{(-1)^k\Gam(z)}{\Gam(z+k)} \intzi D^k(e^{-t\eta(r)}) r^{z+k-1} dr \\
&= \f{(-1)^k\Gam(z)}{\Gam(z+k)} M_t^k(z)
\end{align*}
if each of the boundary terms vanishes.

To see that the boundary terms vanish, we use the condition \eqref{S4-1:E2}  to obtain 
\[|rD(e^{-t\eta(r)})| = |trD\eta(r)| e^{-t\eta(r)} \le Atr^\alpha e^{-t\eta(r)}\]
and
\[|r^2D^2(e^{-t\eta(r)})| \le |tr^2D^2\eta(r)| e^{-t\eta(r)} + |t^2r^2D \eta(r)|^2 e^{-t\eta(r)} \le (Atr^\alpha+(Atr^\alpha)^2) e^{-t\eta(r)}.\]
Inductively, we can obtain for $1\le m\le k$
\begin{equation}\label{S4-1:E4}
|r^mD^m(e^{-t\eta(r)})| \le (Atr^\alpha+\dots+(Atr^\alpha)^m) e^{-t\eta(r)}.
\end{equation}
Consequently we have for $0\le m\le k$
\[|r^{z+m} D^m(e^{-t\eta(r)})| \to 0\]
as $r\to0$.
Using the estimate \eqref{S4-1:E4} and the condition \eqref{S4-1:E1} we also have
\[|r^{z+m} D^m(e^{-t\eta(r)})| \lesssim r^{\Re(z)+\alpha m} e^{-t\eta(r)}  \to 0\]
as $r\to\infty$, and therefore Claim 1 is proved.

\textbf{[Step 2]}
If we show that $M_t(z)$ is integrable, then by Mellin's inversion theorem \eqref{S2-2:E2} we get for $c>0$
\begin{equation}\label{S4-1:E5}
e^{-t\eta(r)} = \frac{(-1)^k}{2\pi i} \int_{(c)} \f{\Gam(z)}{\Gam(z+k)} M_t^k(z) r^{-z} dz.
\end{equation}

Using the estimate \eqref{S4-1:E4} we have
\begin{align*}
|M_t^k(z)|
&\le \intzi |D^k(e^{-t\eta(r)}) r^{z+k-1}| dr \\
&\le kAt \int_{Atr^\alpha<1} r^{\Re(z)-1+\alpha} e^{-t\eta(r)} dr \\
&\quad + k(At)^k \int_{Atr^\alpha>1} r^{\Re(z)-1+\alpha k} e^{-t\eta(r)} dr.
\end{align*}
The integrals on the right converge absolutely because of the condition \eqref{S4-1:E1}.
We notice that for $\Re(z)>-\alpha$ the function $M_t^k(z)$ is bounded and holomorphic on the bounded strip.
Therefore for $k\ge2$
\begin{equation}\label{S4-1:E6}
|M_t(z)| = \Abs{\f{(-1)^k\Gam(z)}{\Gam(z+k)} M_t^k(z)} \lesssim \f{1}{1+|\Im(z)|^k}
\end{equation}
and hence $M_t(z)$ is integrable.

\textbf{[Step 3]}
Using \eqref{S2-1:E1} and \eqref{S4-1:E5}, we have
\begin{align*}
K^\beta(t,x) 
&= (2\pi)^{-d/2} |x|^{-d/2+1} \intzi J_{d/2-1}(|x|r) r^{d/2+\beta} e^{-t\eta(r)} dr \\
&= \f{(-1)^k}{(2\pi)^{d/2} |x|^{d/2-1}} \frac{1}{2\pi i} \intzi \!\!\! \int_{(c)} \f{\Gam(z)}{\Gam(z+k)} M_t^k(z) J_{d/2-1}(|x|r) r^{-z+d/2+\beta} dz dr.
\end{align*}
We choose $(d+1)/2+\beta<c<d+\beta$ so that the integrand
\[\f{\Gam(z)}{\Gam(z+k)} M_t^k(z) J_{d/2-1}(|x|r) r^{-z+d/2+\beta}\]
is integrable from the estimates \eqref{S4-1:E6} and 
\[|J_\nu(r)| \lesssim (r^\nu \wedge r^{-1/2}).\]
So we can apply Fubini's theorem to obtain
\[K^\beta(t,x) = \f{(-1)^k}{(2\pi)^{d/2} |x|^{d/2-1}} \frac{1}{2\pi i} \int_{(c)} \f{\Gam(z)}{\Gam(z+k)} M_t^k(z) \intzi J_{d/2-1}(|x|r) r^{-z+d/2+\beta} dr dz.\]
From the identity \eqref{S2-4:E3}
\begin{align*}
\intzi J_{d/2-1}(|x|r) r^{-z+d/2+\beta} dr
&= |x|^{z-d/2-\beta-1} \intzi r^{-d/2+1} J_{d/2-1}(r) r^{-z+d+\beta-1} dr \\
&= |x|^{z-d/2-\beta-1} \f{2^{-z+d/2+\beta} \Gam(\f{-z+d+\beta}{2})}{\Gam(\f{z-\beta}{2})}.
\end{align*}
This completes the proof. \qed
\end{proof}

\begin{remark}
The assumption \eqref{S4-1:E2} is analogous to that of standard symbol class in the theory of pseudo-differential operators.
\end{remark}

\begin{remark}
The meromorphic continuation of the Mellin transform using resolution of singularities goes back to Bernstein-Gelfand \cite{MR0247457}.
\end{remark}

We note that if $\beta$ is not an even integer, the decay of the kernel is dominated by $\beta$.
Because the intensity of the singularity (in the sense of smoothness) at the origin is important.
We already see this phenomenon from Theorem \ref{S3-2:T2} and Proposition \ref{sec3-thm3}.
This reflects a version of localization principle in Harmonic analysis, that is, if the Fourier transform is regular in the frequency domain, then the function is localized in space domain in the sense of decreasing speed.
This shows that the fractional derivatives are generalized from the Laplacian operator, but they have different features compared with the traditional derivatives.
We consider the case that $\beta$ is an even integer at the end of this section.
We only present the first order approximation for convenience.

\begin{theorem}\label{sec4-thm1}
Suppose that $\eta$ satisfies the same conditions in the previous proposition.
If $\beta \neq 0,2,4,\dots$, then for each $t>0$
\[K^\beta(t,x) \sim \f{2^\beta\Gam(\f{d+\beta}{2})}{\pi^{d/2}\Gam(\f{-\beta}{2})} \f{1}{|x|^{d+\beta}}\]
as $|x|\to\infty$.
\end{theorem}

\begin{proof}
The function $M_t^k(z)$ is holomorphic for $\Re(z)>-\alpha$.
Since $\Gam((-z+d+\beta)/2)$ has no pole for $\Re(z)<d+\beta$ and $\beta$ is not an even integer, the integrand
\begin{equation}\label{S4-1:E10}
\f{\Gam(z)\Gam(\f{-z+d+\beta}{2})2^{-z+\beta}}{\Gam(z+k)\Gam(\f{z-\beta}{2})} M_t^k(z) |x|^z
\end{equation}
in the equation \eqref{S4-1:E3} is meromorphic in $\Re(z)>-\alpha$ and has a simple pole at $z=0$.
We calculate the residue at $z=0$ as follows:
\begin{align*}
\f{\Gam(\f{d+\beta}{2})2^\beta}{\Gam(k)\Gam(\f{-\beta}{2})} M_t^k(0)
&= \f{\Gam(\f{d+\beta}{2})2^\beta}{\Gam(k)\Gam(\f{-\beta}{2})} \intzi D^k(e^{-t\eta(r)}) r^{k-1} dr \\
&= \f{\Gam(\f{d+\beta}{2})2^\beta}{\Gam(\f{-\beta}{2})} (-1)^{k-1} \intzi D(e^{-t\eta(r)}) dr \\
&= \f{\Gam(\f{d+\beta}{2})2^\beta}{\Gam(\f{-\beta}{2})} (-1)^k 
\end{align*}
by integrating by parts and using the Fundamental Theorem of Calculus.
Now we set
\[\tilde{c}=\min\set{\alpha,1}/2\]
and define the contour $\gam_R$ consists of a rectangle with vertices $c-iR, c+iR, -\tilde{c}+iR, -\tilde{c}-iR$ with the counterclockwise orientation.
Since the only pole of the integrand \eqref{S4-1:E10} inside the rectangle $\gam_R$ is at $z=0$, Cauchy's Residue theorem gives
\[\frac{1}{2\pi i} \int_{\gam_R} \f{\Gam(z)\Gam(\f{-z+d+\beta}{2})2^{-z+\beta}}{\Gam(z+k)\Gam(\f{z-\beta}{2})} M_t^k(z) |x|^z dz = \f{(-1)^k\Gam(\f{d+\beta}{2})2^\beta}{\Gam(\f{-\beta}{2})}.\]
Using Stirling's formula and Cauchy's theorem, we can shift the vertical line integration left to obtain
\[K^\beta(t,x) = \f{2^\beta\Gam(\f{d+\beta}{2})}{\pi^{d/2}\Gam(\f{-\beta}{2})} \f{1}{|x|^{d+\beta}} + R(t,x),\]
where
\[R(t,x) = \f{(-1)^k 2^\beta}{\pi^{d/2} |x|^{d+\beta}} \frac{1}{2\pi i} \int_{(-\tilde{c})} \f{\Gam(z)\Gam(\f{-z+d+\beta}{2})2^{-z+\beta}}{\Gam(z+k)\Gam(\f{z-\beta}{2})} M_t^k(z) |x|^z dz.\]
Using Stirling's formula, we have the estimate
\[\Abs{\f{\Gam(z)\Gam(\f{-z+d+\beta}{2})2^{-z+\beta}}{\Gam(z+k)\Gam(\f{z-\beta}{2})} M_t^k(z) |x|^z} \lesssim |\Im(z)|^{-\Re(z)+d/2+\beta-k}.\]
and therefore $|x|^{d+\beta}R(t,x)$ tends to 0 as $|x|$ tends to infinity.
\qed
\end{proof}

\begin{remark}
To clarify the result of Theorem \ref{sec4-thm1}, we wish to add a comment.
The function $K^\beta(t,x)$ depends on $t$, but the first order asymptotic is independent of $t$ as $|x|$ tends to infinity, whenever $\beta$ is not an even integer.
The reason is that if $\beta$ is not an even integer, then the integrand \eqref{S4-1:E10} has trivial simple pole at the origin due to $\Gam(z)$, which is \emph{never cancelled}.
This shows that these L\'evy symbols have little effect on the ultimate behavior of the fractional derivatives of transition functions.
The key point of the proof is that after inserting the Mellin inversion formula \eqref{S4-1:E5} and changing the order of integration, the remaining integrand \eqref{S4-1:E10} has meromorphic continuation.
\end{remark}

\begin{remark}
If $1<\alpha<2$, we choose $\tilde{c}=(\alpha+1)/2$ so that the integrand \eqref{S4-1:E10} has two simple poles at $z=0$ and $z=-1$.
Then the residue at $z=-1$ produces the second order approximation.
Actually, one can go further, if the symbol $\eta$ satisfies the regularity condition for all $k\in\bN$, then we can obtain an asymptotic expansion of the kernel.
\end{remark}

\subsection{The Relaxation of the Condition of the Symbol}
\label{S4-2}
The main contribution to $K^\beta(t,x)$ comes from the integration near $r=0$.
We guess that the assumption \eqref{S4-1:E2} could be weakened away from the origin.
This idea is made precise in the following proposition.

Let $\psi(r)$ be a smooth function defined in $\bR$, with the properties that $\psi(r)=0$ for $r\le1$, and $\psi(r)=1$ for $r\ge2$.
Let us set for $\beta\ge0$
\[\cE^\beta(t,x) = \intoi J_{d/2-1}(|x|r) \psi(r) r^{d/2+\beta} e^{-t\eta(r)} dr.\]
We shall show that this is an error term.

\begin{proposition}\label{sec4-thm2}
Suppose that $\eta(r)\in C^N((0,\infty))$ satisfies the integrability condition
\begin{equation}\label{S4-2:E1}
\intoi r^{d/2-1/2+\beta} e^{-t\eta(r)} dr < \infty.
\end{equation}
Then for fixed $t>0$
\[|\cE^\beta(t,x)| \lesssim |x|^{-N-1/2}\]
as $|x|\to\infty$.
\end{proposition}

\begin{proof}
Fix $t>0$ and denote
\[g(r)=\psi(r) r^{d/2+\beta} e^{-t\eta(r)}.\]
Using \eqref{S2-4:E2}, we have
\begin{align*}
&\intzi J_{\nu-1}(|x|r) g(r) dr \\
&= \intzi (|x|r)^\nu J_{\nu-1}(|x|r) (|x|r)^{-\nu} g(r) dr \\
&= \intzi \f{1}{k} D[(|x|r)^\nu J_{\nu}(|x|r)] (|x|r)^{-\nu} g(r) dr \\
&= \lim_{r\to\infty} \f{1}{k} J_{\nu}(|x|r) g(r) - \f{1}{k} \intzi (|x|r)^\nu J_{\nu}(|x|r) D[(|x|r)^{-\nu} g(r)] dr \\
&= \f{1}{k} \intzi J_{\nu}(|x|r) r^{-1}(\nu I-rD) g(r) dr.
\end{align*}
since $\lim_{r\to\infty} J_{\nu}(|x|r) g(r)=0$.
If we denote
\[L_\nu = r^{-1}(\nu I-rD),\]
then we can write
\[\intzi J_{\nu-1}(|x|r) g(r) dr = \f{1}{|x|} \intzi J_{\nu}(|x|r) L_\nu g(r) dr.\]
Carrying out the repeated integration by parts gives
\[\cE^\beta(t,x) = \f{1}{|x|^N} \intzi J_{\nu+N-1}(|x|r) L_{\nu+N-1} L_{\nu+N-2} \cdots L_\nu  g(r) dr.\]
Using the decay of Bessel functions \eqref{S2-4:E1} and the integrability condition \eqref{S4-2:E1}, we have
\begin{align*}
|\cE^\beta(t,x)|
&\lesssim \f{1}{|x|^{N+1/2}} \intoi r^{-1/2} |L_{\nu+N-1} L_{\nu+N-2} \cdots L_\nu  g(r)| dr \\
&\lesssim \f{1}{|x|^{N+1/2}} \intoi r^{-1/2} \sum_{n=0}^N r^{n-N}|D^n g(r)| dr \\
&\lesssim \f{1}{|x|^{N+1/2}} \sum_{n=0}^N \intoi r^{-1/2+n-N} \sum_{k=0}^n |D^{n-k} r^{d/2+\beta}| |D^k e^{-t\eta(r)}| dr \\
&\lesssim \f{1}{|x|^{N+1/2}} \sum_{n=0}^N \sum_{k=0}^n  \intoi r^{-1/2-N+d/2+\beta+k} |D^k e^{-t\eta(r)}| dr \\
&\lesssim \f{1}{|x|^{N+1/2}}
\end{align*}
and the result follows.
\qed
\end{proof}

We have the same main term under the weak assumptions compare to those of the previous theorem.
This shows that the local behavior of the symbol $\eta$ near $r=0$ is important.

\begin{corollary}\label{sec4-thm3}
Suppose that $\eta(r)$ is a real-valued continuous function on $[0,\infty)$ and that there exist $0<\alpha<2$, $k>(d+3)/2+\beta$, and $M>0$ such that
\begin{equation}\label{S4-2:E2}
\sup_{1\le m\le k} \sup_{0<r<1} r^{-\alpha+m} |D^m \eta(r)| < \infty,
\end{equation}
\begin{equation}\label{S4-2:E3}
\sup_{0\le m\le k} \sup_{1<r<\infty} |D^m \eta(r)| r^{-M} < \infty,
\end{equation}
and
\begin{equation}\label{S4-2:E4}
\liminf_{r\to\infty} \f{\eta(r)}{\log r} = \infty.
\end{equation}
If $\beta \neq 0,2,4,\dots$, then for each $t>0$
\[K^\beta(t,x) \sim \f{2^\beta\Gam(\f{d+\beta}{2})}{\pi^{d/2}\Gam(\f{-\beta}{2})} \f{1}{|x|^{d+\beta}}\]
as $|x|\to\infty$.
\end{corollary}

\begin{proof}
By essentially the same way of in the proof of Theorem \ref{sec4-thm1}, we can prove the following Mellin inversion formula for $c>0$
\[e^{-t\eta(r)} = \frac{(-1)^k}{2\pi i} \int_{(c)} \f{\Gam(z)}{\Gam(z+k)} M_t^k(z) r^{-z} dz,\]
where
\[M_t^k(z) = \intzi u^{z+k-1} \left(\f{d}{du}\right)^k e^{-t\eta(u)} du.\]
The remaining proof is the same and so it is omitted.
\qed
\end{proof}

\subsection{Perturbation of the Symbol}
\label{S4-3}
We consider the perturbation of the symbol corresponding to the symmetric $\alpha$-stable process that has the form $\eta(r)=r^{\alpha}+\eta_1(r)$ under some conditions of $\eta_1(r)$.

\begin{theorem}\label{sec4-thm4}
Let $0<\alpha<2$.
Suppose that $\eta(r)=r^{\alpha}+\eta_1(r)$ and that $\eta_1(r)$ is a real-valued continuous function on $[0,\infty)$ that satisfies for some
$\delta>\alpha$, $k>(d+3)/2+\beta$, and $M>0$ such that
\begin{equation}\label{S4-3:E1}
\sup_{1\le m\le k} \sup_{0<r<1} r^{-\de+m} |D^m \eta_1(r)| < \infty,
\end{equation}
\begin{equation}\label{S4-3:E2}
\sup_{0\le m\le k} \sup_{1<r<\infty} |D^m \eta_1(r)| r^{-M} < \infty,
\end{equation}
and
\begin{equation}\label{S4-3:E3}
\liminf_{r\to\infty} \f{\eta_1(r)}{\log r} = \infty.
\end{equation}
If $\beta=0,2,4,\dots$, then for each $t>0$
\[K^\beta(t,x) \sim -\frac{2^{\beta+\alpha}\Gam(\f{d+\beta+\alpha}{2})}{\pi^{d/2}\Gam(-\f{\beta+\alpha}{2})} \frac{t}{|x|^{d+\beta+\alpha}}\]
as $|x|\to\infty$.
\end{theorem}

\begin{proof}
We already know that
\begin{equation}\label{S4-3:E4}
K^\beta(t,x) = \f{1}{(2\pi)^{d/2} |x|^{d+\beta}} \intzi J_{d/2-1}(r) e^{-t(r/|x|)^\alpha} e^{-t\eta_1(r/|x|)} r^{d/2+\beta} dr.
\end{equation}
Since $\Gam(z)$ is the Mellin transform of $e^{-r}$, we have for $c_0>0$
\begin{equation}\label{S4-3:E5}
e^{-t(r/|x|)^\alpha} = \frac{1}{2\pi i\alpha}\int_ {(c_0)}
\Gam(\f{s}{\alpha}) (t^{1/\alpha}r/|x|)^{-s} ds
\end{equation}
by Mellin's inversion formula \eqref{S2-2:E2}.
We also have for $c>0$
\begin{equation}\label{S4-3:E6}
e^{-t\eta_1(r/|x|)} = \frac{(-1)^k}{2\pi i} \int_{(c)} \f{\Gam(z)}{\Gam(z+k)} M_t^k(z) (r/|x|)^{-z} dz,
\end{equation}
where
\[M_t^k(z) = \intzi u^{z+k-1} \left(\f{d}{du}\right)^k e^{-t\eta_1(u)} du.\]
Putting \eqref{S4-3:E5} and \eqref{S4-3:E6} together with \eqref{S4-3:E4}, we obtain
\begin{align*}
K^\beta(t,x) &= \f{(-1)^{k+1} }{\alpha(2\pi)^{d/2+2}|x|^{d+\beta}}
\intzi \!\!\! \int_{(c_0)} \int_{(c)} \f{\Gam(z)\Gam(\f{s}{\alpha})}{\Gam(z+k)} M_t^k(z)  t^{-s/\alpha}  |x|^{s+z} \\
&\quad \times J_{d/2-1}(r) r^{-s-z+d/2+\beta} dz ds dr.
\end{align*}
From the identity \eqref{S2-4:E3}  we have
\[\intzi J_{d/2-1}(r) r^{-s-z+d/2+\beta} dr = \f{2^{-s-z+d/2+\beta}\Gam(\f{-s-z+d+\beta}{2})}{\Gam(\f{s+z-\beta}{2})}.\]
We choose $(d+1)/2+\beta<c+c_0<d+\beta$ so that we can apply Fubini's theorem and the identity \eqref{S2-4:E3} to obtain
\begin{align*}
p^\beta(t,x) &= \f{(-1)^{k+1} }{\alpha(2\pi)^{d/2+2}|x|^{d+\beta}}  \int_{(c_0)}
\int_{(c)} \f{\Gam(z)\Gam(\f{s}{\alpha})}{\Gam(z+k)} M_t^k(z) t^{-s/\alpha} |x|^{s+z} \\
&\quad \times \intzi J_{d/2-1}(r) r^{-s+z+d/2+\beta} dr dz ds \\
&= \f{(-1)^{k+1} }{\alpha(2\pi)^{d/2+2} |x|^{d+\beta}}
\int_{(c_0)} \int_{(c)} \f{\Gam(z)\Gam(\f{s}{\alpha})}{\Gam(z+k)} M_t^k(z) t^{-s/\alpha} |x|^{s+z} \\
&\quad \times \f{2^{-s-z+d/2+\beta}\Gam(\f{-s-z+d+\beta}{2})}{\Gam(\f{s+z-\beta}{2})} dz ds \\
&= \f{(-1)^{k+1}2^{d/2+\beta} }{\alpha(2\pi)^{d/2+2}|x|^{d+\beta}} \int_{(c_0)}\int_{(c)}
\f{\Gam(s/\alpha)\Gam(z)\Gam(\f{-s-z+d+\beta}{2})}{\Gam(z+k)\Gam(\f{s+z-\beta}{2})} \\
&\quad \times M_t^k(z) t^{-s/\alpha} (|x|/2)^{s+z}dz ds.
\end{align*}

Now, we define the contour $\gam_R$ consists of a
rectangle with vertices $c-iR, c+iR, -\de/2+iR, -\de/2-iR$ with the
counterclockwise orientation.
The integrand
\[\f{\Gam(z)\Gam(\f{-s-z+d+\beta}{2})}{\Gam(z+k)\Gam(\f{s+z-\beta}{2})} M_t^k(z) (|x|/2)^{z}\]
has only one simple pole at $z=0$ in the rectangle bounded by
$\gam_R$. By Cauchy's Residue Theorem,
\[\frac{1}{2\pi i} \int_{\gam_R} \f{\Gam(z)\Gam(\f{-s-z+d+\beta}{2})}{\Gam(z+k)\Gam(\f{s+z-\beta}{2})} M_t^k(z) (|x|/2)^{z} dz = \f
{\Gam(\f{-s+d+\beta}{2})}{\Gam(k)\Gam(\f{s-\beta}{2})} M_t^k(0).\] We
have $M_t^k(0) = (-1)^{k} \Gam(k)e^{-t\eta_1(0)}$.
As $R\to\infty$, the two integrals over the line parallel to $x$-axis go to zero by Stirling's formula. Thus, we can shift the line
integral left to obtain
\[\frac{1}{2\pi i} \int_{(c)} \f{\Gam(z)\Gam(\f{-s-z+d+\beta}{2})}{\Gam(z+k)\Gam(\f{s+z-\beta}{2})} M_t^k(z) (|x|/2)^{s+z} dz = \f {(-1)^k\Gam(\f{-s+d+\beta}{2})}{\Gam(\f{s-\beta}{2})}e^{-t\eta_1(0)} + R_{t,x}(s),\]
where
\[R_{t,x}(s)=\frac{1}{2\pi i} \int_{(-\frac{\alpha+\delta}{2})} \f{\Gam(z)\Gam(\f{-s-z+d+\beta}{2})}{\Gam(z+k)\Gam(\f{s+z-\beta}{2})} M_t^k(z) (|x|/2)^{z} dz.\]
Hence we have
\begin{align*}
K^\beta(t,x)
&= \f{-2^{d/2+\beta}ie^{-t\eta_1(0)}}{\alpha(2\pi)^{d/2+1} |x|^{d+\beta}}
\int_{(c_0)}\f{\Gam(s/\alpha) \Gam(\f{-s+d+\beta}{2})}{ \Gam(\f{s-\beta}{2})} t^{-s/\alpha}(|x|/2)^s ds \\
&\quad + \f{(-1)^{k+1}2^{d/2+\beta}i}{\alpha(2\pi)^{d/2+1}|x|^{d+\beta}} \int_{(c_0)} \Gam(s/\alpha) R_{t,x}(s) t^{-s/\alpha}|x|^{s}   ds\\
&=-\frac{2^{\beta+\alpha}\Gam(\f{d+\beta+\alpha}{2})}{\pi^{d/2}\Gam(-\f{\beta+\alpha}{2})} \frac{te^{-t\eta_1(0)}}{|x|^{d+\beta+\alpha}}\\
&\quad + \f{-2^{d/2+\beta}ie^{-t\eta_1(0)} }{\alpha(2\pi)^{d/2+1}
|x|^{d+\beta}}\int_{(-\frac{3}{2}\alpha)}\f {\Gam(s/\alpha)
\Gam(\f{-s+d+\beta}{2})}{ \Gam(\f{s-\beta}{2})}
 t^{-s/\alpha}(|x|/2)^sds \\
&\quad + \f{(-1)^{k+1}2^{d/2+\beta}i}{\alpha(2\pi)^{d/2+1}|x|^{d+\beta}} \int_{(\frac{\delta-\alpha}{4})} \Gam(s/\alpha)
R_{t,x}(s) t^{-s/\alpha}|x|^{s} ds.
\end{align*}
For $k>(d+3)/2+\beta$, the last two integrals converge absolutely by Stirling's formula and they are bounded by
\[|x|^{-d-\beta-\min\set{\frac{3\alpha+\delta}{4}, \frac{3}{2}\alpha}},\]
as $|x|\to\infty$.
Since $\min\set{\frac{3\alpha+\delta}{4}, \frac{3}{2}\alpha}>\alpha$, the theorem is proved.
\qed
\end{proof}

\section*{Acknowledgement}

The authors thank Zhen-Qing Chen for his helpful comments and for giving related references during the summer school in Probability 2012 in Republic of Korea.

\def\cprime{$'$}

\end{document}